\documentclass[12pt,reqno]{amsart}
\usepackage{amsmath} 
\usepackage{amssymb}

\def\squarebox#1{\hbox to #1{\hfill\vbox to #1{\vfill}}}

\newcommand{\w}{\omega}

\newcommand{\Z}{{\mathbb{Z}}}
\newcommand{\C}{{\mathbb{C}}}
\newcommand{\F}{{\mathcal F}}
\newcommand{\CC}{{C_c^\infty(\rr)}}

\newcommand{\lw}{L_{\|W\|}^p(\rr)}

\newcommand{\W}{W(E)}

\newcommand{\Fv} {{\bf \F}}
\newcommand{\iy}{+\infty}

\newcommand{\R}{{\mathbb R}}
\newcommand{\RR}{{\mathbb R}}
\newcommand{\rr}{{\mathbb R}^+}
\newcommand{\N}{{\mathbb N}}
\newcommand{\T}{{\bf T}}

\newcommand{\dr}{C_c^\infty(\rr)}
\newcommand{\supp}{\operatorname{supp}}

\newcommand{\Ss}{{\bf {S}}}
\newcommand{\uv}{M_{u,v}}
\newcommand{\puv}{\nu_{a,u,v}}
\newcommand{\pt}{\Nu_{a}}

\newcommand{\Cc}{{C_c^{\infty}({\mathbb R}^+)}}

\newcommand{\CH}{{C_c^\infty({\rr})}\otimes H}
\newcommand{\EE}{\overline{E}}
\newcommand{\WE}{W(\overline{E})}
\newcommand{\PP}{{\bf P^+}}
\newcommand{\Nu}{{\mathcal{V}}}

\newcommand{\Sa}{{\bf S}_a}
\newcommand{\Sma}{{\bf S}_{-a}}

\renewcommand{\Im}{\mathop{\rm Im}\nolimits}

\theoremstyle{plain}

\newtheorem{thm}{Theorem}

\newtheorem{lem}{Lemma}

\newtheorem{prop}{Proposition}

\newtheorem{rem}{Remark}
\newtheorem{deff}{Definition}

\theoremstyle{definition}

\numberwithin{equation}{section}

\title{ Wiener-Hopf operators on spaces of functions on $\R^+$ with values in a Hilbert space}
\author[ V. Petkova]{ Violeta Petkova}
\begin{document}\address{Violeta Petkova\\
Universit\'e Paul S\'ebatier, \\
UFR: MIG\\
 Laboratoire Emile Picard \\
 118 route de Narbonne  \\
31062 Toulouse Cedex 4, France.\\}
\email{petkova@math.ups-tlse.fr}

\maketitle
\vspace{0.5cm}
{\bf Abstract.} A Wiener-Hopf operator on a Banach space of functions on $\rr$ is a bounded operator $T$ such that $P^+ S_{-a}TS_{a}=T$, $a\geq 0$, where $S_a$ is the operator of translation by $a$. We obtain a representation theorem for the Wiener-Hopf operators on a large class of functions on $\rr$ with values in a separable Hilbert space. \\ 


\section{Introduction}
This paper deals with Wiener-Hopf operators on Banach spaces of functions on $\R^+$ with values in a separable Hilbert space $H$.  Let $E$ be a Banach space of functions on $\rr$ such that $E\subset L_{loc}^1(\rr)$.
 For $a\geq 0$, define the operator
$$S_a:E\longrightarrow L_{loc}^1(\rr),$$
 by the formula
$(S_a f)(x)=f(x-a)$, for almost every $x\in [a, +\infty[$ and $(S_af)(x)=0$, for $x\in [0,a[$. 
For $a\geq 0$, introduce
$$S_{-a}:E\longrightarrow L_{loc}^1(\rr),$$
 defined by the formula
$(S_{-a}f)(x)=f(x+a),$ for almost every $x\in\rr.$ Notice that $S_{-a}S_a=I$ but $S_aS_{-a}\neq I$.
>From now, we suppose that $S_{a}E\subset E$ and $S_{-a}E\subset E$, $\forall a\in \rr$. 
The Wiener-Hopf operators on $E$ are the bounded operators 
$$T: E\longrightarrow E$$
 satisfying 
$$S_{-a}TS_{a}=T,\:\forall a\in \rr.$$
Denote by $P^+$ the operator
$$P^+: L^1_{loc}(\R)\longrightarrow L_{loc}^1(\rr)$$ 
defined by 
$$(P^+f)(x)=f(x),\:a.e.\: {\rm on} \:\rr.$$
The Wiener-Hopf operators which appear in theory of the signal and in control theory have been studied in a lot of papers. The problem we deal here is the existence of a symbol for operators of this type. 
It is well-known that if $T$ is a Wiener-Hopf operator on $L^2(\rr)$ there exists $h\in L^\infty(\R)$ such that 
$$Tf=P^+\F^{-1}(h\hat{f}),\:\forall f\in L^2(\rr).$$
Here $\F$ denotes the usual Fourier transformation from $L^2(\R)$ into $L^2(\R)$. 
The function $h$ is called the symbol of $T$. 
Despite of the extensive literature related to Wiener-Hopf operators, there are not analogous representation theorem for Wiener-Hopf operators on general Banach spaces of functions even if the functions are with values in $\C$. Here we develop a theory of the existence of a $L^\infty$ symbol for every Wiener-Hopf operator in a very large class of spaces of functions on $\rr$ with values in a separable Hilbert space. Moreover, we obtain a caracterisation of $spec(S_1)\cap (spec(S_{-1}))^{-1}$. 
The determination of the spectrum of a translation operator is an open question in general spaces of functions on $\rr$ and it plays an important role in the scattering theory. 
We are motivated by the results of \cite{V2} proving the existence of a symbol for every Wiener-Hopf operator on a weighted space $L_\w^2(\R^+)$ (see Example 1 for the definition). On the other hand, the methods exposed in \cite{V5} and \cite{G} show that the existence of the symbol of a multiplier (a bounded operator commuting with the translations) on spaces of scalar functions on $\R$ implies an analogous result for the multipliers on a space of functions on $\R$ with values in an Hilbert space. The arguments in \cite{V5} and \cite{G} have been based on the link between the scalar and the vector-valued cases. 
 However the results concerning the symbol of a multiplier do not imply analogous results about Wiener-Hopf operator in the general case. It is well-known that for every Wiener-Hopf operator $T$ on $L^2(\R^+)$, there exists a multiplier $M$ on $L^2(\R)$ such that $P^+M=T$. Unfortunaly, a such result is not known even for Wiener-Hopf operators on a weighed space $L_\w^2(\R^+)$. Despite some progress (see \cite{G}, \cite{V5}) in the study of the symbol of a multiplier on a space of functions on $\R$ with values in a Hilbert space, the analogous problem for Wiener-Hopf operators has been very few considered. Moreover, even in the case of the weighted spaces of functions on $\rr$ with values in a Hilbert space the existence of the symbol of a Wiener-Hopf operator was an open problem still now. 
First in Section 2, we improve the results of \cite{V2} concerning the existence of the symbol of a Wiener-Hopf operator on $L_\w^2(\rr)$ replacing $L_\w^2(\rr)$ by a general Banach space of functions on $\rr$ satisfying only three natural hypothesis given below. Next following the methods of \cite{V5} and \cite{G} and using the results of Section 2, we obtain the existence of the symbol of a Wiener-Hopf operator on a very large class of spaces of functions on $\rr$ with values in a separable Hilbert space. In Section 4 we explain how the setup considered here can by extended in several directions.\\

Let $E$ be a Banach space of functions on $\rr$ with values in $\C$ satisfying the following three hypothesis. \\
(H1) We have $\Cc\subset E\subset L_{loc}^1(\rr)$, the inclusions are continuous and $\Cc$ is dense in $E$.\\
(H2) For every $x\in \RR$, $S_xE=E$ and $\sup_{x\in K}\|S_x\|<+\infty,$ for every compact $K$ of $\R$.\\
(H3) For all $a\in \R$, the operator 
$\Gamma_{a}$ defined by 
$$(\Gamma_{a} f)(x)=e^{iax}f(x),\:a.e.\:,\forall f\in E$$
is bounded on $E$ and
$$\sup_{a\in \R}\|\Gamma_a\|< +\infty.$$

Notice that (H3) is trivial, if we have $\|f\|=\Big \|\:|f|\:\Big\|$ in E.
Let $C_K^\infty(\rr)$ be the space of $C^\infty$ functions with a compact support included in $K$. 
 For simplicity, we will write $S$ instead of $S_1$.
Since the norm of $f$ given by $\sup_{a\in \R} \|\Gamma_a f\|$ is equivalent to the norm of $E$, we will assume from now that $\Gamma_a$ is an isometry for every $a\in\R$. 
Denote by $\rho(A)$ the spectral radius of a bounded operator $A$. 
Set
$$I_E=[-\ln \rho(S_{-1}), \ln \rho(S)]$$
and 
$$U_E=\Big\{ z\in \C,\:\Im z \in I_E\Big\}.$$
For $f\in E$, denote by $(f)_a$ the function defined by $(f)_a(x)=e^{ax}f(x)$, a.e. on $\R^+$.
 In Section 2 we obtain the following result which generalizes Theorem 1 in \cite{V2}.
\begin{thm}
Let $T\in \W$.\\
$1$) For every $a \in I_E$ we have $(Tf)_a \in L^2(\RR^+)$, for $\:f \in \Cc.$\\
$2$) For every $a \in I_E$ there exists a function $\nu_a \in L^{\infty}(\RR)$ such that 
$$(Tf)_a=P^+{\mathcal F}^{-1}(\nu_a \widehat {(f)_a}), \:for\: f \in \Cc$$
and  we have $\|\nu_a\|_\infty\leq C\|T\|,$  where $C$ is a constant dependent only on $E$.\\
$3$) Moreover, if  $\overset{\circ}{I_E} \neq \emptyset$ $(i.e.\:\frac{1}{\rho(S_{-1})}<\rho(S))$, there exists a function $\nu \in {\mathcal H}^\infty(\overset{\circ}{U_E}  )$ such that for every $a \in\overset{\circ}{I_E}$ we have
$$\nu(x+ia)=\nu_a(x),\:\:almost \:everywhere\: on\: \R.$$
\end{thm}

\begin{deff}
If $ \overset{\circ}{I_E} \neq \emptyset$, $\nu$ is called the symbol of $T$, and if $I_E=\{a\}$, then $\nu_a$ is the symbol of $T$.
\end{deff}

Using Theorem 1, we also obtain the following spectral result. 
\begin{thm}
 We have $$spec(S)\cap \Big(spec(S_{-1})\Big)^{-1}=\Big\{ z\in \C,\:\frac{1}{\rho(S_{-1})}\leq |z|\leq \rho(S)\Big\}.$$
\end{thm}

This result is new even in the case of the spaces $L^2_{\omega}(\R^+)$. In particular, we conclude that if $\rho(S) > \frac{1}{\rho(S_{-1})}$ the spectrum of $S$ contains a disk. The proof of Theorem 2 is based on the existence of a symbol for every Wiener-Hopf operator and the construction of suitable cut-off function $f \in C_c^{\infty}(\R^+)$. This application was one of the motivations to search a symbol of a Wiener-Hopf operator. Moreover, we extent below the same result for operators with values in a Hilbert space (see Theorem 4).\\

The main result of this paper is an analogous result for Wiener-Hopf operators on spaces of functions on $\R^+$ with values in a separable Hilbert space. Denote by $<u,v>$ the scalar product of $u,\:v\in H$ and let $\|u\|_H$ be the norm of $u\in H$. 
Denote by $L_{loc}^1(\rr,H)$ the space of functions 
$$F:\rr\longrightarrow H$$ such that 
$$\Big(\rr\ni x\longrightarrow \|F(x)\|_H\Big)\in L_{loc}^1(\rr).$$
Let $\mathcal{L}(H)$ be the space of bounded operators on $H$.
Introduce the vector space $\CH$ generated by $fu$ for $f\in \dr$ and $u\in H$. Denote by $C_0(\R^+,H)$ the Banach space of all norm continuous functions
$$\Phi:\rr \longrightarrow H$$
such that for every $\epsilon>0$, there exists a compact set $K_\epsilon$ such that 
$$\|\Phi(x)\|_H=0,\:\forall x\in \R^+\setminus K_\epsilon.$$
Let $E$ be a Banach space of functions on $\rr$ with values in $\C$ satisfying (H1), (H2) and (H3). 
Denote by $\EE$ the Banach space of functions 
$$F:\rr \longrightarrow H$$
such that
$$\Big(\rr\ni x\longrightarrow \|F(x)\|_H\Big)\in E.$$
We will see in Section 3 that $C_c^\infty(\R^+)\otimes H$ is dense in $\EE$. 
For illustration, we give below some examples.\\

{\bf Example 1.} Let $E=L_\w^p(\R^+)$,
where $\w$ is a weight on $\rr$ and $p\in [1,+\infty[.$
We recall that $w$ is a weight on $\rr$ if $\w$ is a non-negative measurable function on $\rr$ such that for all $y\in\rr$,
$$0<\sup_{x\in \rr}\frac{\w(x+y)}{\w(x)}<+\infty$$
and
$$0<\sup_{x\in \rr}\frac{\w(x)}{\w(x+y)}<+\infty.$$ 
The space $L_\w^p(\R^+)$ is the set of measurable functions $f$ from $\rr$ into $\C$ such that 
$$\int_{\rr} |f(x)|^p\w(x)^pdx <+\infty ,$$ 
equipped with the norm 
$$\|f\|_{\w,p}=\Big( \int_{\rr} |f(x)|^p\w(x)^pdx \Big)^{\frac{1}{p}}.$$
It is easy to see that $L_\w^p(\rr)$ satisfies the hypothesis (H1), (H2) and (H3). For the study of the Wiener-Hopf operators on $L_\w^2(\rr)$ the reader may consult \cite{V2}. The space $\overline{E}$ associated to $L_\w^p(\rr)$ is the space usually denoted by $L_\w^p(\rr,H)$ of functions 
$$F:\rr\longrightarrow H$$
such that 
$$\int_{\rr} \|F(x)\|_H^p\:\w(x)^pdx<+\infty.$$

{\bf Example 2.}
Let $A$ be a real-valued continuous function on $[0, +\infty[$, such that $A(0)=0$ and let $\frac{A(y)}{y}$ be non-decreasing for $y>0$. Let
$L_A(\rr)$ be the set of all complex-valued, measurable functions on $\rr$ such that 
$$\int_{\rr} A\Big( \frac{|f(x)|}{t}\Big)dx<+\infty,$$
for some positive number $t$ and let
$$\|f\|_A=\inf\Big\{t>0\:|\:\int_{\rr} A\Big( \frac{|f(x)|}{t}\Big)dx\leq 1\Big\},$$
for $f\in L_A(\rr)$. Then $L_A(\rr)$ is a Banach space called a Birnbaum-Orlicz space (see \cite{B}). It is easy to check that $L_A(\rr)$ satisfies (H1), (H2) and (H3).
If $E=L_A(\rr)$, the associated space $\overline{E}$ is the set $L_A(\rr, H)$ of measurable functions 
$$F:\rr\longrightarrow H$$
such that for some $t>0$, we have
$$\int_{\rr}  A\Big( \frac{\|F(x)\|_H}{t}\Big)dx<+\infty.$$

{\bf Example 3.}
Let $A$ be a function satisfying the properties described in Example 2. Let $\w$ be a weight on $\rr$. Define $L_{A,\w}(\rr)$ as the space of measurable functions on $\rr$ such that 
$$\int_{\rr} A\Big( \frac{|f(x)|}{t}\Big)\w(x)dx<+\infty,$$
for some positive number $t$ and let
$$\|f\|_{A,\w}=\inf\Big\{t>0\:|\:\int_{\rr} A\Big( \frac{|f(x)|}{t}\Big)\w(x)dx\leq 1\Big\}$$
for $f\in L_{A,\w}(\rr)$. Then $L_{A,\w}(\rr)$ is a Banach space called a weighted Orlicz space. It is easy to check that $L_{A,\w}(\rr)$ satisfies (H1), (H2) and (H3). If $E=L_{A,\w}(\rr)$, the associated space $\overline{E}$ is the set $L_{A,\w}(\rr, H)$ of measurable functions 
$$F:\rr\longrightarrow H$$
such that for some $t>0$,
$$\int_{\rr}  A\Big( \frac{\|F(x)\|_H}{t}\Big)\w(x)dx<+\infty.$$

For $a>0$, we define the operators 
$$\Sa:\EE\longrightarrow \EE$$
and 
$$\Sma:\EE\longrightarrow \EE$$
by 
$$(\Sa F)(x)=F(x-a),\:a.e. \:on \:[a,+\infty[,$$
$$(\Sa F)(x)=0,\:\forall x\in [0, a[,$$
$$(\Sma F)(x)=F(x+a),\:a.e.\:on\:\rr.$$
For simplicity, we will write $\Ss$ instead of $\Ss_1$. For $F\in \overline{E}$, we denote by $\|F\|_H$ the function
$$\|F\|_H:\rr\ni x\longrightarrow \|F(x)\|_H\in \C.$$
For fixed $a\in \R$, we see that for $F\in \EE$, $F\neq 0$, we have
$$\frac{\|\Sa F\|}{\|F \|}=\frac {\:\|S_a(\|F\|_H)\:\|}{\|\:\|F\|_H\:\|}\leq \|S_a\|.$$
We conclude that $\Sa$ is bounded and $\| \Sa\|\leq \|S_a\|$. If $\|f\|=\|\:|f|\:\|$, for every $f\in E$, obviously we get $\| \Sa\|= \|S_a\|$.
Introduce the operator 
$$\PP: L^1_{loc}(\R, H)\longrightarrow L^1_{loc}(\rr, H)$$
 defined by the formula
$$(\PP F)(x)=F(x), \:a.e. \:{\rm on} \:\rr.$$

\begin{deff}
We call a Wiener-Hopf operator on $\EE$ every bounded operator $\T$ on $\EE$ such that
$$\T  \Phi=\Sma  \T \Sa  \Phi,\:\forall a>0,\:\forall \Phi\in \EE.$$
Denote by $\WE$ the set of the Wiener-Hopf operators on $\EE$.
\end{deff}
The main result of this paper is the following. 
\begin{thm}
Let $E$ be a Banach space satyisfing (H1), (H2) and (H3). Let $\T\in \WE$.\\
$1$) We have $(\T \Phi)_a\in L^2(\rr,H),\:\forall \Phi\in \CH,\:\forall a\in I_E.$\\
$2$) There exists $\Nu_a\in L^\infty(\RR,\mathcal{L}(H))$ such that
$$(\T \Phi)_a=\PP\F^{-1}(\Nu_a (.)[\widehat{(\Phi)_a}(.)]),\:\forall a\in I_E,\:\forall \Phi\in \CH.$$
Moreover, ${{\rm{ess}}}\:\sup_{x\in \R} \|\Nu_a(x)\|\leq C\|\T\|,$ where $C$ is a constant dependent only on $E$.\\
$3$) If $\overset{\circ}{U_E}\neq \emptyset$, set
$$\Nu(x+ia)=\Nu_a(x),\:\forall a \in\overset{\circ}{I_E},\:for\:almost \:every\:x \in \R.$$
Then for $u$, $v\in H$, the function
$$z\longrightarrow <u, \Nu(z)[v]>$$
is in ${\mathcal H}^\infty(\overset{\circ}{U_E})$
and $\sup_{z\in\overset{\circ}{U_E}  }\|\Nu(z)\|\leq C \|\T\|.$
\end{thm}
\begin{rem}
We will see later that $\rho(\Ss)=\rho(S)$, $\rho(\Ss_{-1})=\rho(S_{-1)}$ and 
$$I_E=[-\ln \rho(S_{-1}), \ln \rho(S)]=[-\ln \rho(\Ss_{-1}), \ln \rho(\Ss)].$$
\end{rem}
We also obtain the following.
\begin{thm}
 We have
$$spec(\Ss)\cap \Big(spec(\Ss_{-1})\Big)^{-1}=\Big\{z\in \C,\:\frac{1}{\rho(\Ss_{-1})}\leq |z|\leq \rho(\Ss)\Big\}.$$
\end{thm}

The spectral cracterisation in Therorem 4 has not been known until now even in particular cases when $E$ is a weighted $L^p$ space with a simple weight.  
\vspace{0.5cm}

\section{Wiener-Hopf operators on Banach spaces of scalar functions on $\R^+$}
\vspace{0.5cm}

In this section, we prove Theorem 1. We follow  the arguments of \cite{V2} in our more general case. For the reader convenience we give the details of the steps which need some modifications.
First, we show that every Wiener-Hopf operator is associated to a distribution. Denote by $C_0^\infty(\rr)$ the space of functions of $C^\infty(\RR)$ with support in $]0, \iy[$. Set 
$$H^1(\RR)=\{f \in L^2(\RR)\:|\:f^{\prime} \in  L^2(\RR)\},$$
 the derivative of $f \in L^2(\RR)$ being computed in the sense of distributions. 
\begin{lem}
If $T \in \W$ and $f \in C_K^\infty(\rr)$, then $(Tf)^{\prime}=T(f^{\prime}).$ 
\end{lem}
{\bf Proof.}
Let  $f \in C_K^\infty(\rr)$ and let $(h_n)_{n \geq 0}\subset \RR$ be a sequence converging to 0. Since 
$\frac{S_{h_n}f-f}{h_n}$ converges to $f^\prime$ with respect to the topology of $C_K^\infty(\rr)$, $\frac{S_{h_n}f-f}{h_n}$ converges to $f^\prime$ with respect to the topology of $E$. Then we have 
$$\int_{\rr} T(f^\prime)(x)\phi(x)dx=\int_{\rr} T\Big(\lim_{n\to +\infty}\frac{S_{h_n}f-f}{h_n}\Big)(x)\phi(x)dx$$
$$=\int_{\rr} \lim_{n\to +\infty}\Big(\frac{S_{h_n}Tf-Tf}{h_n}\Big)(x)\phi(x)dx=
\int_{\rr} (Tf)(x)\phi^\prime(x)dx,\:\forall \phi \in \dr.$$
Consequentelly,  $T(f^{\prime})=(Tf)^{\prime}$ in the sense of distributions.  
$\Box$
\begin{prop}
If $T$ is a Wiener-Hopf operator, then there exists a distribution $\mu_T$ of order 1 such that 
$$Tf=P^+(\mu_T*f),$$
for $f \in \Cc$.
\end{prop}

The proof of Proposition 1 follows the arguments of that of Theorem 2 in \cite{V2} and we omit it. We just give the definition of $\mu_T$. We have $$<\mu_T,f>=\lim_{x \to \iy} (TS_x{\tilde f})(x),$$
for $f\in C_c^\infty(\R)$, where $\tilde{f}$ is the function defined by ${\tilde f}(x)=f(-x),$ for $f\in C_c^\infty(\R)$, $x \in \RR$.

\begin{deff}
If $\phi \in C_c^\infty(\R)$, we denote by $T_\phi$ the Wiener-Hopf operator such that 
$$T_\phi f=P^+ (\phi*f),\:\forall f\in C_c^\infty(\rr).$$
\end{deff}

\begin{prop}
 If $T\in \W$, then there exists a sequence $(\phi_n)_{n \in\N}\subset \dr$ such that 
 $$\lim_{n \to \iy}\|T_{\phi_n}f-Tf\|=0, \forall f \in E$$
 and
 $$\|T_{\phi_n}\|\leq C\|T\|\:,\forall n \in \N,$$
where $C$ is a constant depending only on $E$.
 \end{prop}

{\bf Proof.} The proof follows the idea of the proof of Theorem 3 in \cite{V2}, but here we must work with Bochner integrals and this leads to some difficulties. For the conviviance of the reader we give the details. 
Let $T \in \W$ and set ${\mathcal T}(t)=\Gamma_{t}\circ T \circ \Gamma_{-t}, \forall t \in  \RR$.  For $a>0$, and $f \in E$ we have
$$ (S_{-a}{\mathcal T}(t)S_af)(x)=({\mathcal T}(t)S_af)(x+a)$$
$$=e^{it(x+a)}\Big(T(f(s-a)e^{-its})\Big)(x+a)$$
$$=e^{itx}\Big(S_{-a}T\Big(f(s-a)e^{-it(s-a)}\Big)\Big)(x)$$
$$=e^{itx}(S_{-a}TS_a(\Gamma_{-t}f))(x)=({\mathcal T}(t)f)(x), \:a.e.$$
This shows that ${\mathcal T}(t)\in \W$. Moreover, we have $\|{\mathcal T}(t)\|=\|T\|$, for $t \in \RR$ and ${\mathcal T}(0)=T$. The application ${\mathcal T}$ is continuous from $\RR$ into $\W$. For $n \in \N$, $ \eta \in \RR$, $x \in \RR$, set 
$$g_n(\eta):=\Bigl(1-\Big|\frac{\eta}{n}\Big|\Bigr)\chi_{[-n,n]}(\eta)$$
 and 
$$\gamma_n (x)=\frac {1-cos(nx)}{ \pi x^2n}.$$
 We have $ \widehat {\gamma _n }(\eta)=g_n(\eta),\:\:\forall \eta \in \RR,\:\:\forall n \in \N.$ Clearly, 
$\|\gamma_n \|_{L^1}=1$ for all $n$ and 
$$\lim_{n \to +\infty}\int_{|x| \geq a}\gamma_n(x)dx=0,\:\forall a>0.$$
Set $Y_n:=({\mathcal T} * \gamma_n)(0).$ Then for $f \in E$ we obtain 
$$\lim_{n \to \iy}\|Y_nf-Tf\|=0.$$
We claim that for $f\in C_K^\infty(\rr)$, we have
\begin{equation}\label{eq:cl}
(Y_nf)(y)=\int_{\RR} ({\mathcal T}(x)f)(y)\gamma_n(-x)dx,\:\forall y \in \rr
\end{equation}
>From Lemma 1, we know that for fixed $x\in \RR$ the function
$$\rr\ni y\longrightarrow ({\mathcal T}(x)f)(y)$$
is $C^\infty$. Let $K_0$ be a compact subset of $\rr$ and let $\psi\in C_{K_0}^\infty(\rr)$. We see that
$$|\psi(y) ({\mathcal T}(x)f)(y)|=|\psi(y)(\mu_T*\Gamma_{-x}(f))(y)|=
|\psi(y)<\mu_{T,z},f(z-y)e^{-ix(z-y)}>|$$
$$\leq C(\mu)\|\psi\|_\infty(\|S_y\Gamma_{-x}f\|_\infty+\|(S_y\Gamma_{-x}f)^\prime\|_\infty)$$
$$\leq C(\mu)\|\psi\|_\infty(\|f\|_\infty+\|f^\prime \|_\infty),\:\forall y\in K_0.$$
Consequently,
$$\int_{\RR} \|\psi {\mathcal T}(x)f\|_\infty \gamma_n(-x) dx <+\infty$$
 and hence the integral
$$\int_{\RR} \psi ({\mathcal T}(x)f)\gamma_n(-x)dx$$ 
is a well-defined Bochner integral with values in $C_{K_0}^\infty(\rr)$. The map 
$$C_{K_0}^\infty(\rr)\ni g \longrightarrow g(x)\in \C,$$
is a continuous linear form for every $x\in \rr$. Since Bochner integrals commute with continuous linear forms (see \cite{HP}) we have
$$\psi(y) (Y_nf)(y)=\psi(y)\int_{\RR} ({\mathcal T}(x)f)(y)\gamma_n(-x)dx,\:\forall y\in \rr$$ 
and the claim (\ref{eq:cl}) is proved.\\
It is clear that
$$\|Y_nf\|\leq \int_{\RR} \| \mathcal{T}(x)f\|\gamma_n(-x)dx\leq \|\mathcal{T}(x)\| \|f\|,\:\forall f\in E.$$
Since $\|\mathcal{T}(x)\|=\|T\|$, we get $\|Y_n\|\leq \|T\|,\:\forall n\in \N.$ \\

 Now consider the distribution associated to $Y_n$. Let $K$ be  a compact subset of $\RR$ and let $z_K \geq 1$ be such that $K \subset ]-\infty, z_K[$. 
 Choose $g \in \dr$ such  that $g$ is positive, $\supp g \subset [z_K-1, z_K+1]$ and $g(z_K)=1.$ 
For $f \in C_K^{\infty}(\RR)$, we have $gT(S_{z_K}({\tilde f}g_n))\in H^1(\RR)$ and it follows from Sobolev's lemma (see \cite{R}) that
$$|(TS_{z_K}({\tilde f}g_n))(z_K)|=|g(z_K)(TS_{z_K}({\tilde f}g_n))(z_K)|$$
$$\leq C\Bigl( \Bigl(\int_{|y-z_K|\leq 1}g(y)^2|(TS_{z_K}({\tilde f}g_n)(y)|^2dy \Bigr)^{\frac{1}{2}}+
 \Bigl(\int_{|y-z_K|\leq 1}|(g(TS_{z_K}({\tilde f}g_n))^{\prime}(y)|^2dy \Bigr)^{\frac{1}{2}}\Bigr),$$
where $C>0$ is a constant. Taking into account (H1), $T$ may be considered as a bounded operator from $\dr$ into $L_{loc}^1(\rr)$ and we have
$$|(TS_{z_K}({\tilde f}g_n))(z_K)|$$
$$\leq C(K)\|T\|(\|\tilde{f}g_n\|_\infty+\| (\tilde{f}g_n)^\prime\|_\infty)$$
$$\leq {\tilde C}(K)(\|f\|_\infty+\|f^{\prime}\|_\infty),$$
  where $C(K)$ and ${\tilde C}(K)$ are constants depending only on $K$. Therefore
 $$|(TS_{z}({\tilde f}g_n))(z)|\leq {\tilde C}(K)(\|f\|_\infty+\|f^{\prime}\|_\infty), \:\forall z\geq z_K, \:\forall f \in C_K^{\infty}(\RR)$$
and we conclude that $\mu_T g_n$ defined by
  $$<\mu_T g_n, f>=\lim_{z \to \iy}(TS_z({\tilde f}g_n))(z)$$
 is a distribution of order 1. On the other hand,  we have
 $$(Y_nf)(y)=\int_{\RR}({\mathcal T}(-s)f)(y)\gamma_n(s)ds=\int_\RR e^{-isy}(T(\Gamma_{s}f))(y)\gamma_n(s)ds$$
 $$=\int_{\RR}<\mu_{T,x},f(y-x)e^{-isx}>\gamma_n(s)ds=<\mu_{T,x},f(y-x)\int_{\RR} \gamma_n(s)e^{-isx}ds>$$
 $$=<\mu_{T,x},f(y-x)g_n(x)>=(\mu_T g_n*f)(y), \forall y\geq 0, \: \forall f \in \Cc.$$
 Finally, we obtain $$Y_nf=P^+(\mu_T g_n*f), \: \forall f \in \Cc,\:\forall n \in \N .$$
Since $\supp \mu_Tg_n \subset [-n,n]$, it is sufficient to obtain the Proposition 2 for $T\in \W$ such that $\mu_T$ is a distribution with compact support. Without lost of generality we assume that $\mu_T$ is with compact support. 
 Let  $(\theta_n)_{n \in \N}\subset \dr$ be a sequence such that $\supp \theta_n \subset [0, \frac{1}{n}]$, $\theta_n \geq 0$,
$$\lim_{n \to  \iy}\int_{x\geq a} \theta_n (x)dx=0,\:\forall a >0$$
 and $\|\theta_n \|_{L^1}=1,$ for $n \in \N$. For $f \in E$ we have 
$$\lim_{n \to \iy} \|\theta_n *f-f\|=0.$$
 Set 
$$T_nf=T(\theta_n *f), \:\forall f \in E.$$
 We conclude that $(T_n)_{n \in \N}$ converges to $T$ with respect to the strong operator topology and $T_n=T_{\phi_n},$ where $\phi_n=\mu_T*\theta _n\in \dr.$ For $f \in E$, we have
 $$\|T_nf\|=\Big\|P^+\Big(\int_0^{\frac{1}{n}}\theta_n(y)S_y(\mu_T*f)dy\Big)\Big\|$$
$$\leq \Big\|\int_0^{\frac{1}{n}}\theta_n(y)P^+(\mu_T*S_yf)dy\Big\|$$
$$\leq \int_0^{\frac{1}{n}}\theta_n(y)\|T\|\|S_y\| \|f\| dy,\:\forall f\in \dr.$$
Then we obtain 
$$\|T_n\|\leq \Big(\int_0^{\frac{1}{n}}\theta_n(y)\|S_y\|\:dy\Big)\|T\|,\:\forall n\in \N$$
and this completes the proof of the proposition. 
$\Box$\\

We need also the following lemma.


\begin{lem}
For every $\phi\in \dr$, we have
\begin{equation}{\label{eq:p}}
|\hat{\phi}(\alpha)|\leq \|T_\phi\| ,\:\forall \alpha \in U_E.
\end{equation}

\end{lem}
{\bf Proof.} 
We use the fact that for a bounded operator $A$ on $E$, there exists a sequence $(f_{n})_{n\in \N}\subset E$ such that:
\begin{equation}\label{eq:f}
\lim_{n\to +\infty}\|A f_n-\rho(A) f_n\|=0\:{\rm and}\:\|f_n\|=1,\:\forall n\in \N.
\end{equation}
Fix $\lambda=\rho(S)$. Let $(f_{n,1})_{n\in \N}$ be a sequence of $E$ such that 
$$\lim_{n\to +\infty}\|S f_{n,1}-\rho(S) f_{n,1}\|=0$$
and
$$\|f_{n,1}\|=1,\:\forall n\in \N.$$
For $p\in\N^*$, observe that 
$$\lambda^{\frac{1}{p}}=\rho(S_{\frac{1}{p}}).$$
Let $(f_{n,\frac{1}{p}})_{n\in\N}\subset E$ be a sequence such that 
$$\lim_{n\to +\infty}\Big\|S_{\frac{1}{p}} f_{n,\frac{1}{p}}-\rho(S_{\frac{1}{p}}) f_{n,\frac{1}{p}}\Big\|=0$$
and
$$\:\Big\|f_{n,\frac{1}{p}}\Big\|=1,\:\forall n\in\N.$$
Notice that for all $q \in {\Bbb N}^*$, such that $q\leq p$ we have:
$$ \Big\| S_{\frac{1}{q}}\:f_{n,\frac{1}{p!}}\:-\:\lambda^{\frac{1}{q}}\:f_{n,\frac{1}{p!}} \Big\| \:=\:
 \Big\| (S_{\frac{1}{p!}})^{\frac{p!}{q}}\:\:f_{n,\frac{1}{p!}}\:-\:(\lambda^{\frac{1}{p!}})^{\frac{p!}{q}}\:\:f_{n,\frac{1}{p!}} \Big\|$$
$$ \leq \: \Big(\prod_{u \in {\Bbb C},\:u^{\frac{p!}{q}}=1,\:u \neq 1} \: \Big\| S_{\frac{1}{p!}}\:-\:u\lambda^{\frac{1}{p!}} \Big\|\:\: \Big)
\Big\| S_{\frac{1}{p!}}\:f_{n,\frac{1}{p!}}\:-\:\lambda^{\frac{1}{p!}}\:f_{n,\frac{1}{p!}}\Big \|. $$ 
We have
$$ \prod_{ u \in {\Bbb C},\:u^{\frac{p!}{q}}=1,\:u \neq 1}\: \Big\| S_{\frac{1}{p!}}\:-\:u\lambda^{\frac{1}{p!}}\Big\|\leq C,\:$$
where $C$ is a constant independent of $n$ and hence we have  
$$\lim_{n \to+\infty} \Big\| S_{\frac{1}{q}}\:f_{n,\frac{1}{p!}}-\lambda^{\frac{1}{q}}\:f_{n,\frac{1}{p!}}\Big \| \:=\:0.$$
Consequentelly, by a diagonal extraction, we can construct $(f_{n})_{n \in \N}$ such that :
$$\lim_{ n\to+\infty}\Big\| S_{\frac{1}{p}}\:f_{n}\:-\:\lambda^{\frac{1}{p}}\:f_{n}\Big\| \:=\:0,\:\:\forall p\in {{\Bbb N}^*} $$
and
$$\| f_{n}\|\:=\:1, \:\:\forall n\in {\Bbb N}.$$ 
For all $p \in {\Bbb N}^*$ and for all $q \in {\Bbb N}$, we have
$$S_{\frac{q}{p}}\:f_{n}-\lambda^{\frac{ q}{p}}\:f_{n}={\mathcal C}_{ q,p}\:\:(S_{\frac{1}{p}}-\lambda^{\frac{1}{p}}I)\:f_{ n},$$
 where ${\mathcal C}_{ q,p}$ is a linear combination of translations. Then
$$\|S_{\frac{q}{p}}\:f_{n}-\lambda^{\frac{ q}{p}}\:f_{n}\| \leq \|\mathcal{C}_{ q,p}\|\:\:\|S_{\frac{1}{p}}\:f_{ n}-\lambda^{\frac{1}{p}}\:f_{n}\|,\:\:\forall n \in \N$$
and
$$\lim_{n\to +\infty}\|S_{\frac{q}{p}}\:f_n-\lambda^{\frac{q}{p}}f_n\|=0.$$
On the other hand,
$$\|S_{-\frac{q}{p}}\:f_{n}-\lambda^{-\frac{ q}{p}}\:f_{n}\| \leq|\lambda^{-\frac{ q}{p}}| \:\|S_{-\frac{q}{p}}\|\:\:\|\lambda^{\frac{p}{q}}f_{n}-S_{\frac{q}{p}}\:f_{ n}\|,\:\:\forall n \in \N $$
and
$$\lim_{ n\to +\infty}\|S_{-\frac{q}{p}}\:f_{n}-\lambda^{-\frac{q}{p}}\:f_{n}\| =0,\:\forall p\in \N^*,\:\:\forall q\in \N .$$
Since $\Bbb Q$ is dense in $\Bbb R$, we deduce that
$$\lim_{n \to + \infty}\|S_{t}\:f_{n}-\lambda^{t}\:f_{n}\| =0,\:\:\forall t\in {\Bbb R}. $$
Now, fix $\phi\in \dr$. Notice that 
$$\int_\R \phi(x) S_x f_n \: dx$$ is a well-defined Bochner interval on $E$ and 
\begin{equation}\label{eq:B}
T_\phi f_n=\int_\R \phi(x) \:(S_x f_n)\: dx.
\end{equation}
Indeed, let $K$ be a compact subset of $\rr$. We have $T_\phi(C_K^\infty(\R^+))\subset C^\infty_{K+supp(\phi)}(\rr)$ and the restriction of $\int_\R\phi(x)S_xdx$ to $C_K^\infty(\rr)$ can be considered as a Bochner integral on $C_K^\infty(\rr)$ with values in $C_{K+supp(\phi)}(\rr)$. It is clear that for $x\in \rr$, the map
$$f\longrightarrow f(x)$$
is a continuous linear form on $C_{K_0}^\infty(\rr)$, for every compact $K_0$. 
Since Bochner integrals commute with continuous linear forms, we obtain, for $g\in \CC$,
$$(T_\phi g)(x)=(\phi*g)(x)=\int_{\RR} \phi(y)g(x-y)dy=\int_{supp(\phi)}\phi(y)(S_yg)(x)dy$$
$$=\Big(\int_{supp(\phi)}\phi(y)S_y g \Big )(x),\:\forall x\in \rr$$
and the formula (\ref{eq:B}) follows from the density of $\CC$ in $E$. 
Then, for all $n\in \N$, we get
$$\Big|\int_{\R} \phi(x)\lambda^x dx\Big|=\Big \| \Big(\int_{\R} \phi(x)\lambda^x dx\Big) f_n\Big \|$$
$$\leq \Big \|\int_{\R} \phi(x)\lambda^x f_n dx-\int_{\R} \phi(x)S_x f_n dx\Big\|+\Big\| \int_{\R} \phi(x)S_x f_n dx\Big\|$$
$$\leq \int_{\R} |\phi(x)|\|\lambda^x f_n-S_x f_n\| dx +\|T_\phi\|.$$
Taking into account the properties of $(f_n)_{n\in\N}$ and the dominated convergence theorem, it follows that 
$$\lim_{n\to +\infty}  \int_{\R} |\phi(x)|\|\lambda^x f_n-S_x f_n\| dx =0.$$
Denote by $C_r$ the circle of radius $r$ and denote by $\mathcal{D}_r$ the line
$$\mathcal{D}_r=\{z\in \C\:|\: \Im z=r\}.$$
We will write $e^{ia.}$ for the function
$$x\longrightarrow e^{iax}.$$
Since $\|T_\phi\|=\|T_{e^{ia.}\phi}\|,$ for all $a\in \R$,
we obtain 
$$\Big|\int_{\R} \phi(x)\lambda^x dx\Big|\leq \|T_\phi\|,\:\forall \lambda \in C_{\rho(S)}.$$
We conclude that 
$$|\hat{\phi}(\alpha)|\leq \|T_\phi\|,\:\forall \phi \in C_c^\infty(\R),\:\forall \alpha\in \mathcal{D}_{\ln \rho(S)}.$$
Denote by $E^*$ the dual space of $E$ and denote by $\|.\|_*$ the norm of $E^*$. 
For $\lambda\in C_{\frac{1}{\rho(S_{-1})}}$, applying the same methods in $E^*$, we obtain that there exists a sequence $(g_{n})_{n\in\N}\subset E^*$ such that 
$$\lim_{n\to +\infty}\|(S_x)^* g_n-\lambda^{x} g_n\|_*=0$$
and 
$$\|g_n\|_*=1,\:\forall n\in \N.$$
We notice that we have 
\begin{equation}\label{eq:aj}
T_\phi^*=\Big(\int_{\rr}\phi(x) S_x dx\Big)^*=\int_{\rr} \phi(x)(S_x)^*dx,\:\forall \phi\in C_c^\infty(\R),
\end{equation}
see \cite{HP}. Then we obtain  as above that 
$$|\hat{\phi}(\alpha)|\leq \|T_\phi^*\|=\|T_\phi\|,\:\forall \phi \in C_c^\infty(\R),\:\forall \alpha\in \mathcal{D}_{-\ln \rho(S_{-1})}.$$
>From the Phragmen-Lindel\"of theorem, it follows that 
$$|\hat{\phi}(\alpha)|\leq \|T_\phi\|,\:\forall \phi \in C_c^\infty(\R),\:\forall \alpha \in U_E.$$
$\Box$\\

The proof of Theorem 1 follows from Proposition 2 and Lemma 2 exactly in the same way as in the proof of Theorem 1 in \cite{V2} and we omit the details.\\

In the proof of Theorem 2, we need the following technical lemma. 
\begin{lem}
Let $\epsilon > 0, \eta_0 > 0$ and $V = \{\xi \in \R^+: |\eta_0 - \xi| \leq \delta\} \subset \R^+$ be fixed. Let $C_0 > 0$ be a fixed constant.
 For $t_0 > 0$ sufficiently large there exists a function $f \in C_c^{\infty}(\R^+)$ with the properties:
\begin{equation} \label{eq:1.1}
 \int_{\R \setminus V} |\hat{f}(\xi)| d \xi \leq \epsilon /C_0.
\end{equation}
\begin{equation} \label{eq:1.2}
 \int_{\R} |\hat{f}(\xi)| d \xi \leq 2\sqrt{2\pi}.
\end{equation}
\begin{equation}\label{eq:1.3}
|f(t_0)| = 1.
\end{equation}
\end{lem}
{\bf Proof.} Introduce the function $g$ with Fourier transform
$$\hat{g}(\xi) = \frac{1}{a} e^{-\frac{(\xi - \eta_0)^2}{2 a^2}} e^{-i t_0 \xi},$$
where $a > 0$ will be taken small enough below. We have
$$\int_{\R \setminus V} |\hat{g}(\xi)| d \xi = \frac{1}{a} \int_{|\xi - \eta_0| \geq \delta} e^{-\frac{(\xi - \eta_0)^2}{2 a^2}} d \xi$$
$$\leq e^{-\frac{\delta^2}{4 a^2}} \frac{1}{a} \int_{\R} e^{-\frac{(\xi - \eta_0)^2}{4 a^2}}d \xi \leq \frac{\epsilon}{2C_0}$$
for $ a > 0$ small enough. We fix $a > 0$ with this property. Obviously,
$$\int_{\R} |\hat{g}(\xi)| d\xi = \int_{\R} e^{-\mu^2/2}d\mu = \sqrt{2 \pi}.$$ 
On the other hand,
$$g(t) = \frac{1}{2 \pi a} \int_{\R} e^{-\frac{(\xi - \eta_0)^2}{2 a^2}} e^{i(t - t_0)\xi} d\xi$$
$$ = \frac{1}{2\pi} e^{i(t - t_0)\eta_0} \int_{\R} e^{-\mu^2 / 2} e^{i(t - t_0) a \mu} d \mu = e^{-\frac{ a^2 (t - t_0)^2}{2}} e^{i(t- t_0)\eta_0}$$ 
and $|g(t_0)| = 1.$\\

Now we will take $t_0 > 2$ sufficiently large.  Let $\varphi \in C_c^{\infty}(\R)$ be a fixed function such that
$\varphi (t) = 0$ for $t \leq 1/2$ and for $t \geq 2t_0 - 1/2$ and let $\varphi(t) = 1$ for $1 \leq t \leq 2t_0 - 1, \: 0 \leq \varphi \leq 1.$ Introduce the function $f = \varphi g \in C_c^{\infty}(\R^+).$ The property (\ref{eq:1.3}) is trivial. We will show that (\ref{eq:1.1}) is satisfied for $t_0 > 0$ large enough depending on the choice of $a > 0$. The proof of (\ref{eq:1.2}) is similar and easier.\\

The function $F= (\varphi - 1)g$ has a small Fourier transform. Moreover, given $\epsilon > 0$ we can take $t_0 > 0$ large enough in order to have
\begin{equation} \label{eq:1.4}
|(1 + \xi^2) \hat{F}(\xi)| \leq \frac{\epsilon }{2 \pi C_0}.
\end{equation}
Indeed, for $\xi^2 \hat{F}(\xi)$ we use an integration by parts with respect to $t$ using the fact that
$$\xi^2 e^{-it\xi} = -\partial_{t}^2 \Bigl(e^{-it \xi}\Bigr).$$
On the support of $(\varphi - 1)$ we have $|t - t_0| > t_0 - 1.$ Thus after the integration by parts in the integral 
$\int_\R e^{-it\xi}(1+\xi^2)F(t)dt$
we are going to estimate an integral
$$\int_{|t - t_0| \geq t_0 -1} e^{-\frac{a^2(t - t_0)^2}{2}}|P(t)| dt$$
with $P$ a polynomial of degree not greater than 2. 

To get (\ref{eq:1.4}), remark that this integral is bounded by 
$$C \Bigl[\int_{-\infty}^{1 - t_0} y^2 e^{-a^2 y^2 /2} dy + \int_{t_0 - 1}^{\infty} y^2 e^{-a^2 y^2 / 2} dy \Bigr]$$
and taking $t_0 > 0$ sufficiently large we arrange (\ref{eq:1.4}). Next we obtain
$$\int_{\R \setminus V} |\hat{f}(\xi)| d \xi \leq \int_{\R \setminus V_{\epsilon}} |\hat{g}(\xi)| d\xi + \int_{\R \setminus V_{\epsilon}} |\hat{F}(\xi)| d \xi$$
$$\leq \frac{\epsilon}{2 C_0} + \frac{\epsilon}{2 \pi C_0} \int_{\R} (1 + \xi^2)^{-1} d\xi \leq \frac{\epsilon}{C_0}.$$
The proof of the lemma is complete. $\Box$\\

{\bf Proof of Theorem 2.}
 First, we show that
\begin{equation}\label{eq:S1}\Big\{ z\in \C,\:\frac{1}{\rho(S_{-1})}\leq |z|\leq \rho(S)\Big\}\subset spec(S).
\end{equation}
Fix $\lambda\notin spec(S)$. Then the operator $(S-\lambda I)^{-1}$ is a Wiener-Hopf operator and following 2) of Theorem 1, we get
$$(S-\lambda I)^{-1} (f)_a=P^+ \F^{-1}(\nu_a\widehat{(f)_a)},\:\forall a\in I_E,\:\forall f\in \dr,$$
where $\nu_a \in L^\infty(\R)$. 
Replacing, $f$ by $(S-\lambda I)g$, we obtain
$$(g)_a=P^+\F^{-1}\Big(\nu_a\F\Big(((S-\lambda I)g)_a\Big)\Big), \:\forall g\in \dr.$$
Denote by $e^{a+i.}$ the function
$$x\longrightarrow e^{a+ix}.$$
It is easy to see that 
$$\F((Sg)_a)(t)=e^{a-it}\F((g)_a)(t),\:\forall a\in I_E,\:\forall t\in \R,\forall g\in \dr.
$$
Consequently,
\begin{equation} \label{eq:for}(g)_a(t)=\F^{-1}[(e^{a-i.}-\lambda) \nu_a\widehat{(g)_a})](t),\:\forall a\in I_E,\:\forall t\in \rr,\forall g\in \dr.
\end{equation}
We have 
\begin{equation}\label{eq:inq}
\|\F^{-1}[ (e^{a-i.}-\lambda)   \nu_a\widehat{(g)_a}]\|_\infty\leq\||(e^{a-i.}-\lambda) |\nu_a\widehat{(g)_a}\|_{L^1(\R)},\:\forall a\in I_E,\:\forall g\in \dr.
\end{equation}
Now, suppose that $|\lambda|=e^b$, for some $b\in I_E$. Choose a small $\epsilon\in ]0, 1[$. It is easy to find an interval ${V_\epsilon}\subset \rr$ such that 
$$|e^{b-it}-\lambda|\leq \frac{\epsilon}{2\|\nu_b\|_\infty},\:\forall t\in V_\epsilon.$$
 Taking into account Lemma 3, we can choose $g\in \dr$ satisfying the following three conditions:\\
1) $\int_{\R\setminus{V_\epsilon}} |\widehat{(g)_b}(t)|dt\leq \frac{\epsilon}{4e^{b}\|\nu_b\|_\infty}$\\
2) 
$\int_{V_\epsilon} |\widehat{(g)_b}(t)|dt\leq 1$\\
3) There exists $t_0\in \rr$, such that $|(g)_b(t_0)|\geq \epsilon$.\\
Taking into account that (\ref{eq:for}) and (\ref{eq:inq}) hold for $g\in \dr$,we get
$$|(g)_b(t_0)|\leq \int_{V_\epsilon} |e^{b-it}-\lambda|\|\nu_b\|_\infty|\widehat{(g)_b}(t)|dt +  \int_{\R\setminus{V_\epsilon}} |e^{b-it}-\lambda|\|\nu_b\|_\infty|\widehat{(g)_b}(t)|dt\leq \epsilon.$$
Hence we obtain a contradiction so 
$|\lambda|\neq e^{a}$, $\forall a\in I_E$ and (\ref{eq:S1}) follows.


We will prove now that
\begin{equation}\label{eq:S-1}
\Big\{z\in \C,\:\frac{1}{\rho(S)}\leq |z|\leq \rho(S_{-1})\Big\}\subset spec(S_{-1}).
\end{equation} 
Let $\lambda \notin spec(S_{-1})$. Then $(S_{-1}-\lambda I)^{-1}\in \W. $ Indeed, for all $x\in \rr$, we observe that
$$S_{-x} (S_{-1}-\lambda I)^{-1} S_x$$
$$=(S_{-1}-\lambda I)^{-1}(S_{-1}-\lambda I)S_{-x}(S_{-1}-\lambda I)^{-1} S_x$$
$$=(S_{-1}-\lambda I)^{-1}S_{-x}(S_{-1}-\lambda I)(S_{-1}-\lambda I)^{-1}S_x$$
$$=(S_{-1}-\lambda I)^{-1}.$$
Hence, for all $g\in \dr$ and for each $a\in I_E$, we have 
$$((S_{-1}-\lambda I)^{-1}g)_a=P^+\F^{-1}(h_a\widehat{(g)_a}),$$
for some $h_a\in L^\infty(\R)$
and 
$$(f)_a=P^+\F^{-1} \Big(h_a \F(((S_{-1}-\lambda I)f)_a)\Big), \:\forall f\in \dr.$$
Then 
$$\F\Big(((S_{-1}-\lambda I)f)_a\Big)(t)=(e^{it-a}-\lambda) \widehat{(f)_a}(t),\: a.e.\:{\rm on} \:\rr,$$
if we suppose that $supp(f)\subset [1,\infty[$. Repeating the argument of the proof of (\ref{eq:S1}), we get a contradiction if $|\lambda|=e^{-a}$, for some $a\in I_E$. We conclude that
$$\Big\{z\in \C,\:\frac{1}{\rho(S)}\leq |z|\leq \rho(S_{-1})\Big\}\subset spec(S_{-1}).$$
It follows that, if $z\in \C$ is such that $\frac{1}{\rho(S_{-1})}\leq |z|\leq \rho(S)$ then $\frac{1}{z}\in spec(S_{-1})$ and
we deduce that  
$$\Big\{z\in \C,\:\frac{1}{\rho(S_{-1})}\leq |z|\leq \rho(S)\Big\}\subset spec(S)\cap \Big(spec(S_{-1})\Big)^{-1}.$$
>From the definition of the spectral radius we get immediately that 
$$spec(S)\cap \Big(spec(S_{-1})\Big)^{-1}\subset \Big\{z\in \C,\:\frac{1}{\rho(S_{-1})}\leq |z|\leq \rho(S)\Big\}$$ and the proof of Theorem 2 is complete. 
$\Box$

\begin{prop}
If $\phi\in C_c^\infty(\R)$, then 
$$\widehat{\phi}(U_E)\subset spec(T_\phi).$$
\end{prop}
{\bf Proof.}
Fix $\lambda\notin spec(T_\phi)$. Then $(T_\phi-\lambda I)^{-1}$ is a Wiener-Hopf operator and we obtain as above 
$$(g)_a=P^+\F^{-1}(\nu_a[\widehat{(\phi)_a}-\lambda]\widehat{(g)_a}),\:\forall g\in \dr,\:\forall a\in I_E,$$ where $\nu_a\in L^\infty(\R)$.
Choosing a suitable $g\in \dr$, we obtain in the same way as in the proof of (\ref{eq:S1}) a contradiction if 
$$\widehat{(\phi)_a}(t)=\lambda,$$
for some $a\in I_E$ and some $t\in \R$ and the proposition follows imediatelly. 
$\Box$\\

\section{Wiener-Hopf operators on Banach spaces of functions on $\R^+$ with values in a Hilbert space $H$}
\vspace{0.5cm}

Now let $H$ be a separable Hilbert space. Denote by $<u,v>$ the scalar product of $u,\:v\in H$. Let $\|u\|_H$ be the norm of $u\in H$. In this section we prove Theorem 3.
Let $E$ be the Banach space of functions from $\rr$ into $H$ satisfying (H1), (H2) and (H3). 
Let $\EE$ be a Banach space of functions 
$$F:\rr \longrightarrow H$$
such that
$$\Big(\rr\ni x\longrightarrow \|F(x)\|_H\Big)\in E.$$
We have the following two lemmas.
\begin{lem}
The space $\Cc\otimes H$ is dense in $\overline{E}$.
\end{lem}

{\bf Proof.} Let $\Phi\in \EE$. Then there exists a positive sequence $(\phi_n)_{n\in \N}\subset \dr$ such that 
$$\lim_{n\to +\infty} \|\phi_n-\|\Phi(.)\|_H\|_E=0.$$
For almost every $x\in \rr$, set $$\Phi_n(x)=\phi_n(x)\frac{\Phi(x)}{\|\Phi(x)\|},\:\:{\rm if}\:\Phi(x)\neq 0,$$
$$\Phi_n(x)=0,\:\:{\rm if}\:\Phi(x)=0.$$
We have
$$\|\Phi_n-\Phi\|_{\EE}=\|\:\:\|\Phi_n(.)-\Phi(.)\|_H\:\:\|_E=\|\phi_n-\|\Phi(.)\|_H\|_E$$
and it is clear that 
$$\lim_{n\to +\infty}   \|\Phi_n-\Phi\|_{\EE}         =0.$$
Since $\Cc\otimes H$ is dense in $C_0(\rr,H)$, the space $\Cc\otimes H$ is dense in $\overline{E}$.
$\Box$

\begin{lem}
If $\Phi\in \EE$ and $u\in H$, then the function defined by 
$$\rr\ni x\longrightarrow <\Phi(x),u>\in \C$$
is a element of $E$.
\end{lem}

{\bf Proof.} 
Let $\Big(\sum_{n=1}^N \phi_n u_n\Big)_{N\geq0}$ be a sequence in $\CH$ such that $\phi_n\in \dr$, $\forall n\in \N$ and  
$$\lim_{N\to +\infty}\Big\|\sum_{n=1}^N \phi_n u_n-\Phi\Big\|_{\EE}=0.$$ 
Let $u\in H$.
Then we have 
$$\lim_{N\to +\infty}\Big\|<\sum_{n=1}^N \phi_n(.) u_n,u>-<\Phi(.),u>\Big\|_{E}$$
$$\leq \lim_{N\to +\infty}\Big \|\:\:\Big\|\sum_{n=1}^N\phi_n(.)u_n-\Phi(.)\Big\|_H\:\Big\|_E \|u\|=0.$$
Now, it is clear that 
$$x\longrightarrow <\Phi(x),u>\in \C$$
is a element of $E$.\:\:\:\:\:$\Box$\\

In the proof of Theorem 3 we will also use the following lemma.
\begin{lem}
Let $G\in L^2(\R,H)$ and $v\in H$. 
Then we have 
$$\F(<G(.),v>)(x)=<\Fv(G)(x),v>,$$
for almost every $x\in \R$.
\end{lem}
The reader may find the proof of Lemma 6 in \cite{V5}. Next we pass to the proof of our main result. \\

{\bf Proof of Theorem 3.}
Let $\T\in \WE$. Fix $u,v\in H$. Define
$T_{u,v}$ on $E$ by the formula
$$(T_{u,v}f)(x)=<\T(fu)(x),v>,\:\forall f\in E,\:a.e.$$
>From Lemma 5, it follows that $T_{u,v}$ is an operator from $E$ into $E$. It is clear that
$$S_{-x}<\T(S_xf u),v>=<S_{-x}\T(S_xf u),v>=<T(fu),v>,\:\forall x\in \R^+.$$
Then we see that $T_{u,v}\in W(E)$. Following Theorem 1, for $a \in I_E$ there exists a function $\puv \in L^{\infty}(\RR)$ such that 
$$(T_{u,v}f)_a=P^+{\mathcal F}^{-1}(\puv\widehat {(f)_a}), \:\forall\: f \in \Cc.$$
Let $\mathcal{B}$ be an orthonormal basis of $H$ and let $\mathcal{O}$ be the set of finite linear combinations of elements of $\mathcal{B}$. We have 
$$|\puv(x)|\leq C \|T_{u,v}\|,\:\forall x \in \R\backslash N_{u,v},$$
where $N_{u,v}$ is a set of measure zero. Without loss of generality, we can modify $\puv$ on $N=\cup_{(u,v)\in \mathcal{O}\times \mathcal{O}} N_{u,v}$ in order to obtain
$$|\puv(x)|\leq C \|\uv\|\leq C \|\T\| \|u\|\|v\|,\:\forall u,\:v\in \mathcal{O}, a.e.$$
For fixed $x\in \R\backslash N$ we observe that
$$\mathcal{O}\times\mathcal{O}\ni (u,v)\longrightarrow \puv(x)\in \C$$
is a sesquilinear and continuous form on $\mathcal{O}\times \mathcal{O}$ and since $\mathcal{O}$ is dense in $H$, we conclude that there exists an unique map 
$$H\times H \ni (u,v)\longrightarrow \widetilde{\nu}_{a,u,v}(x)\in \C$$
such that 
$$\widetilde{\nu}_{a,u,v}(x)={\nu}_{a,u,v}(x),\:\forall u,v\in \mathcal{O}.$$
Consequently, there exists an unique map 
$$ \Nu_a:\R\longrightarrow \mathcal{L}(H)$$ 
such that
$$<\Nu_a(x)[u],v>=\widetilde{\nu}_{a,u,v}(x),\:\forall u,\:v\in H,\:a.e.$$
It is clear that 
$$\|\Nu_a(x)\|= \sup_{\|u\|=1,\|v\|=1}|<
\Nu_a(x)[u],v>|\leq C \|T\|,\:a.e.$$
 Fix $a\in I_E $ and $f\in \dr$. It is obvious that we have $\widehat{(f)_a}(x)u\in H,\:\forall x\in \R$. 
Next for almost every $x \in \rr$, we obtain
$$\F^{-1}\Big(<\pt(.)[\widehat{(f)_a}(.)u],v>\Big)(x)=\F^{-1}\Big(<\pt(.)[u],v>\:\widehat{(f)_a}(.)\Big)(x) $$
$$=\F^{-1}\Big(\tilde{\nu}_{a,u,v}(.)\widehat{(f)_a}(.)\Big)(x)=(T_{u,v}f)_a(x).$$
Consequently,
\begin{equation} \label{eq:lem}\F^{-1}(<\pt(.)[\widehat{(f)_a}(.)u],v>)(x)=(<\T[fu](.),v>)_a(x),\end{equation}
for almost every $x\in \rr$. 
Now, consider the function  $\Psi_a$ on $\rr$ defined for almost every $x\in \rr$ by the formula 
$$\Psi_a(x)=\pt(x)[\widehat{(f)_a}(x)u]$$ and observe that $\Psi_a\in L^2(\rr,H)$. 
Indeed, we have 
$$\int_{\rr} \| \pt(x)[\widehat{(f)_a}(x)u]\|^2dx $$
$$\leq \int_{\rr} \| \pt(x)\|^2 \|\widehat{(f)_a}(x)u\|^2dx$$
$$\leq C^2 \|\T\|^2 \int_{\rr} |\widehat{(f)_a}(x)|^2\|u\|^2 dx <+\infty.$$
This makes possible to apply Lemma 6, and we get
$$\F^{-1}(<\pt(.)[\widehat{(f)_a}(.)u],v>)(x)=<\Fv^{-1}(\pt(.)[\widehat{(f)_a}(.)u])(x),v>,$$
for almost every $x\in \rr$. It follows from (\ref{eq:lem}) that we have 
$$(\T[fu])_a(x)=\Fv^{-1}(\pt(.)[\widehat{(f)_a}(.)u])(x),$$
for almost every $x\in \rr$ and this yields
$$(\T[fu])_a\in L^2(\rr,H).$$
This completes the proof of 1) and 2). The proof of 3) uses the same argument as the proof of the assertion 3) of Theorem 1. 
$\Box$\\

{\bf Proof of Theorem 4.}
Fix $\alpha\in \C$ and suppose that $\alpha \notin spec(\Ss)$. Then we have 
$(\Ss-\alpha I)^{-1}\in \WE$ and from Theorem 3, we get 
$$((\Ss-\alpha I)^{-1}F)_a=\F^{-1}(\Nu_a(.)[\widehat{(F)_a}(.)]),\:\forall a\in I_E,\:\forall F\in \CH.$$
Replacing $F$ by $(\Ss-\alpha I)G$, we get
$$(G)_a(x)=\F^{-1}(\Nu_a(.)\F[(\Ss-\alpha I) G](.))(x)=\F^{-1}\Big(\Nu_a(.)[(e^{a-i.}-\alpha)\widehat{(G)_a}(.)]\Big)(x).$$
We have
$$\|(G)_a\|_\infty\leq \|\Nu_a(.)[(e^{a-i.}-\alpha)\widehat{(G)_a}(.)]\|_{L^1(\R)},\:\forall a\in I_E.$$
Then if $|\alpha|=e^a$, for some $a\in I_E$ choosing a suitable $G\in \CH$ in the same way as in the proof of Theorem 2, we obtain a contradiction. Hence, 
$$\Big\{ z\in \C,\:\frac{1}{\rho(S_{-1})}\leq |z|\leq \rho(S)\Big\}\subset spec(\Ss).$$ In the same way, we obtain 
$$\Big\{ z\in \C,\:\frac{1}{\rho(S_{-1})}\leq |z|\leq \rho(S)\Big\}\subset (spec(\Ss_{-1}))^{-1}.$$
It follows that 
$$\Big\{ z\in \C,\:\frac{1}{\rho(S_{-1})}\leq |z|\leq \rho(S)\Big\}\subset spec(\Ss)\cap \Big(spec(\Ss_{-1})\Big)^{-1}.$$
Taking into account that 
$$spec(\Ss)\cap \Big(spec(\Ss_{-1})\Big)^{-1}\subset \Big\{z\in \C,\:\frac{1}{\rho(\Ss_{-1})}\leq |z|\leq \rho(\Ss)\Big\}$$
and 
$$\|\Ss\|\leq \|S\|,\:\:\:\| \Ss_{-1}\|\leq \|S_{-1}\|,$$
(see Section 1), we observe that 
$$\rho(\Ss)=\rho(S),\:\:\rho(\Ss_{-1})=\rho(S_{-1}).$$
We deduce that 
$$spec(\Ss)\cap \Big(spec(\Ss_{-1})\Big)^{-1}= \Big\{z\in \C,\:\frac{1}{\rho(\Ss_{-1})}\leq |z|\leq \rho(\Ss)\Big\}$$
and the proof of Theorem 4 is complete. 
$\Box$

\section{Generalizations}
\vspace{0.5cm}
In this section we first deal with the Wiener-Hopf operators in a lager class of Banach spaces of functions on $\rr$ with values in a separable Hilbert space. 
Let $W$ be an operator-valued weight on $\rr$. It means that 
$$W: \rr \longrightarrow \mathcal{L}(H)$$
and $W$ satisfies the property
\begin{equation}\label{eq:prop}0<\sup_{x\in \rr} \frac{\|W(x+y)\|}{\|W(x)\|}<+\infty,\:\forall y\in \rr.
\end{equation}
This implies (see \cite{V1}, \cite{V2}) that for every compact $K$ of $\rr$, we have
$$\sup_{x\in K} \|W(x)\|<+\infty.$$ Notice that if $H$ has a finite dimension, $W$ is given by a matrix. 
We denote by $L_W^p(\rr,H)$ the space of measurable functions $F$ on $\rr$ with values in $H$ such that 
$$\int_{\rr} \|W(x)[F(x)]\|_H^p dx <+\infty. $$
For illustration we give a simple example. \\

{\bf Example.} If $H$ is the space $\R^5$, the operator-valued weight $W$ defined for $x$ by the matrix
$$\left( \begin{array}{ccccc}
  1& e^x &e^{3x}&1&1 \\ 1+x& x&e^x&1&e^{3x}\\ e^x&1&1&x&x+1\\1&1&e^x&e^{2x}&1\\x&x&1+x&e^x&\frac{x^2}{2}
  \end{array} \right) 
  $$
 is such that
 the condition (\ref{eq:prop}) trivially holds.\\

The space $L_W^p(\rr,H)$  is equipped with the norm
$$\Big(\int_{\rr} \|W(x)[F(x)]\|_H^p dx\Big)^\frac{1}{p}. $$
Let $\T$ be a Wiener-Hopf operator on $L_W^p(\rr,H)$. We fix $u$, $v\in H$. 
Notice that for $u\in H$ and $f\in \lw$, we have $fu \in L_W^p(\rr,H)$. Indeed,
$$\int_{\rr} \|W(x)[f(x)u]\|^p dx\leq \int_{\rr} \|W(x)\|^p |f(x)|^p \|u\|^p dx<+\infty.$$
Introduce the operator
$T_{u,v}$ defined on $\lw$ by the formula
$$(T_{u,v}f)(x)=<\T(fu)(x),v>,\:a.e.,\:\forall f\in \lw.$$
It is easy to see that
$$\int_{\rr}\| W(x)\|^p\:|<\T(fu)(x),v>|^p\: dx$$
$$\leq \int_{\rr} \|W(x)\|^p \|\T\|^p \|u\|^p
|f(x)|^p\|v\|dx<+\infty.$$
Consequently, 
$T_{u,v}$ is a Wiener-Hopf operator on $\lw$. Therefore $T_{u,v}$ has a symbol following Theorem 1. Applying the methods exposed in Section 3, we obtain that Theorem 3 holds also if we replace $\EE$ by $L_W^p(\rr,H)$, for $1\leq p <\infty$. Denote by $I_{W}$ (resp. $U_W$) the set $I_E$ (resp. $U_E$) for $E=\lw$. We recall that $I_E$ and $U_E$ are defined in the Introduction. We have the following. 
\begin{thm}
 Let $\T$ by a Wiener-Hopf operator on $L_W^p(\rr,H)$, for $1\leq p <\infty$. \\
$1$) We have $(\T \Phi)_a\in L^2(\rr,H),\:\forall \Phi\in \CH,\:\forall a\in I_{W}.$\\
$2$) There exists $\Nu_a\in L^\infty(\RR,\mathcal{L}(H))$ such that
$$(\T \Phi)_a=\PP\F^{-1}(\Nu_a (.)[\widehat{(\Phi)_a}](.)),\:\forall a\in I_{W},\:\forall \Phi\in \CH.$$
Moreover, ${{\rm{ess}}}\:\sup_{x\in \R} \|\Nu_a(x)\|\leq C\|\T\|.$\\
$3$) If $\overset{\circ}{U_{W}}\neq \emptyset$, 
 set $$\Nu(x+ia)=\Nu_a(x),\:\forall a \in\overset{\circ}{I_{W}},\:for\:almost \:every\:x \in \R$$
We have $\sup_{z\in \overset{\circ}{U_{W}}}\|\Nu(z)\|\leq C \|\T\|$
and for $u$, $v\in H$, the function
$$z\longrightarrow <\Nu(z)u,v>$$
is analytic on $\overset{\circ}{U_{W}}$.\\
\end{thm}

The results of Section 3 and Section 4 hold if we replace $H$ by a separable Banach space $B$ satisfying the following conditions:\\
1) B has a countable basis.\\ 
2) The dual space of $B$ denoted by $B^*$ has a countable basis.\\
For example these conditions are satisfied if $B=l_\w^p(\Z)$, where $\w$ is a weight on $\Z$ and $1\leq p<+\infty$. 
We recall that $\w$ is a weight on $\Z$, if $\w$ is a positive sequence on $\Z$ satisfying
$$0<\sup_{k\in \Z}\frac{\w(k+n)}{\w(k)}< +\infty,\:\forall n\in \Z.$$
It is easy to see that $B^*=l_{\w^*}^q(\Z)$, where $q$ is such that $\frac{1}{p}+\frac{1}{q}=1$.
The weight $\w^*$ is given by the formula 
$$\w^*(n)=\frac{1}{\w(-n)}, \:\forall n\in \Z.$$ 
Denote by $e_n$ the sequence defined by $e_n(k)=0$ if $n\neq k$ and $e_n(n)=1$. 
Considering the family $\{e_n\}_{n\in \Z}$ included in $l^p_\w(\Z)$ and in $l_{\w^*}^q(\Z)$, it is trivial to see that les conditions 1) and 2) are satisfied. \\
Let $B$ be a Banach space satisfying 1) and 2). Let $E$ by a Banach space of functions on $\rr$ satisfying (H1)-(H3). 
Denote by $<\:,\:>_B$ the duality between $B$ and $B^*$. 
Let $\EE$ be the space of functions 
$$F:\rr\longrightarrow B$$
such that
$\|F(.)\|_B\in E$. 
Let $\T$ be a Wiener-Hopf operator on $\EE$. Then using the operators $T_{u,v}$ defined by 
$$(T_{u,v}f)(x)=<\T(fu)(x),v>_B,\: \forall u\in B,\:\forall v\in B^*,\:a.e.$$
and the arguments of the proof of Theorem 3, we obtain an extended version of Theorem 3 in the case of spaces of functions on $\rr$ with values in $B$. 
For example Theorem 3 holds for the Wiener-Hopf operators on spaces of the form $L_{\w_1}^p\Big(\rr, l_{\w_2}^q(\Z)\Big)$, for $1\leq p<\infty$, $1\leq q<\infty$, where $\w_1$ (resp. $\w_2$) is a weight on $\rr$ (resp. $\Z$). The arguments developed in this paper do not hold if we replace $l_{\w_2}^q(\Z)$ by $L_{\w_2}^q(\R)$, for $q\neq 2$. The existence of the symbol of a Wiener-Hopf operator on the family of spaces $L_{\w_1}^p\Big(\rr, L_{\w_2}^q(\R)\Big)$ for $q\neq 2$ is an interesting direction of investigation.
\vspace{0.4cm}



\end{document}